\documentclass[a4paper,11pt]{article}

\usepackage[dvips]{graphicx} 
\usepackage[latin1]{inputenc}
\usepackage{a4wide,epic}
\usepackage{amssymb}         
\usepackage{amsmath}         
\usepackage{amsfonts}
\usepackage[usenames,dvipsnames]{xcolor}
\usepackage{curves}

\begin{document}

\title{Integral control on nonlinear spaces:\\ two extensions} 

\author{Zhifei Zhang\thanks{Key Laboratory of Advanced Control and Optimization for Chemical Processes, Ministry of Education, East China University of Science and Technology, No.130, Meilong Road, Shanghai 200237, China; and Department of Electronics \& Information Systems, Ghent University, Belgium (e-mail: zhifei.zhang@ugent.be).},~~Zhihao Ling\thanks{Key Laboratory of Advanced Control and Optimization for Chemical Processes, ECUST Shanghai.}~~\& Alain Sarlette\thanks{QUANTIC lab, INRIA Paris. 2 rue Simone Iff, 75012 Paris, France; and Data Science Lab, Ghent University. Technologiepark 914, 9052 Zwijnaarde, Belgium. tel: +33 1 80 49 43 59; alain.sarlette@inria.fr}}
\date{\today}

\maketitle

\begin{abstract}
This paper applies the recently developed framework for integral control on nonlinear spaces to two non-standard cases. First, we show that the property of perfect target stabilization in presence of actuation bias holds also if this bias is state dependent. This might not be surprising, but for practical purposes it provides an easy way to robustly cancel nonlinear dynamics of the uncontrolled plant. We specifically illustrate this for robust stabilization of a pendulum at arbitrary angle, a problem posed as non-trivial by some colleagues. Second, as previous work has been restricted to systems with as many control inputs as configuration dimensions, we here provide results for integral control of a non-holonomic system. More precisely, we design robust steering control of a rigid body under velocity bias.
\end{abstract}




\section{Introduction}\label{sec:intro}

Integral control, i.e.~adding a feedback term that is proportional to the time-integral of the deviation between actual and target values, is a basic disturbance rejection tool~\cite{MurrayAstromBook}. It lowers the sensitivity to low-frequency disturbances and in particular for a \emph{constant} disturbance input, it allows the system to stabilize perfectly on the target value. Integral control requires minimal knowledge about the system dynamics --- just enough to tune the integral gain --- and thereby is inherently robust to model uncertainties.

On a nonlinear state space like the circle, the torus or the set of rotations, the definition of integral control must be revised. Indeed, integration is canonically defined only for arguments belonging to vector spaces: the integral (or even the sum) of e.g.~different rotation matrices is not a standard defined thing. In general, when coordinates like the Euler angles are used to describe the position on the nonlinear space, the result depends on the parametrization and can only be valid locally. Whereas when the manifold is embedded in a vector space, i.e.~integrating component-wise the $n \times n$ rotation matrices, the result does not belong to the original space.

Recently, in \cite{MaithripalaJournal,MaithripalaCDC} and \cite{OurSCL}, a canonical definition of integral control has been given for systems evolving \emph{on Lie groups}, to which the above examples belong. It consists in integrating the input command, instead of the output. Indeed the inputs at different times usually belong to a tangent bundle, and thanks to Lie group properties these can be mapped uniquely into a reference vector space --- the Lie algebra --- to define integration in an inherent way. Typically the input command contains something ``like proportional action'', and integrating it gives rise to term precisely equivalent to the standard integral control  for linear systems. The above papers establish the benefits of such integral control on Lie groups for fully actuated systems (number of inputs equals dimension of the Lie group). They show that it perfectly rejects an input bias that is constant (when mapped back to the Lie algebra) and give estimations of the basin of attraction of the target equilibrium, which for topological reasons on manifolds can mostly not be the full state space.

The present paper considers two extensions of this framework with regard to its practical use. \newline
$\bullet$ First, we rigorously establish robustness of this integral control on Lie groups to a state-dependent bias. I.e.~if the disturbance to be rejected is constant in time at every point of the configuration space, but varies continuously from point to point in the configuration space (e.g.~on the circle), then integral control can still perfectly reject it. This is not big news in the traditional linear setting, yet it is worth establishing for the Lie group case as well. Moreover it highlights the ability of integral control to reject the effect of \emph{unknown} nonlinear plant dynamics on the steady state. This complements the control approach by cancellation of open-loop dynamics, which requires perfect plant knowledge.\newline
$\bullet$ Second, we investigate integral control on Lie groups for so-called under-actuated or non-holonomic systems, i.e.~where the number of control inputs is strictly lower than the dimension of the Lie group. We do this on a benchmark model, namely a vehicle under steering control \cite{Krishnaprasad}. We propose a new viewpoint on the constant input bias in this case, and a related adaptation of integral control, for which we establish convergence properties. Although our paper may suggest a framework towards a general solution, those results are still weaker than in the case of full actuation and a full treatment of under-actuated systems remains an open question at this point.\\

The paper is organized as follows. In Section \ref{sec:statedep} we analyze the benefit of integral control for state-dependent bias. In Section \ref{ssec:pendulum} we apply this to the stabilization of the nonlinear pendulum at an arbitrary angle; this is a problem of ``simple, robust'' control suggested to us by prof.~R.Sepulchre \cite{RodTutorial}. In Section \ref{sec:car} we analyze the effect of a bias in translation velocity on the steering controlled vehicle. We observe that a problem appears only if the velocity \emph{direction} is not correctly estimated. In Section \ref{sec:car2} we propose and analyze an appropriate integral controller to mitigate this effect.


\section{State-dependent bias}\label{sec:statedep}

Let us first consider a linear system
$$\tfrac{d}{dt} x = Ax + B u + F v \;\;, \;\;\; y = C x$$
where $F v$ is a constant input disturbance, and controlled with a PID controller $u = -k_P y - k_D \tfrac{d}{dt} y - k_I \int y(t) \, dt $. It is well-known that under appropriate stability conditions, the presence of the integral control term makes the system converge to $y=0$ despite the disturbance $v$ \cite{MurrayAstromBook}. It is maybe less known that this can still hold if $F$ depends on $y$. Indeed, if $F=F_0+ F_1\, y$ we have equivalently
$$\tfrac{d}{dt} x = A'x + B u + F_0 v \;\;, \;\;\; y = C x$$
with $A' = A + F_1 C v$. Thus as long as the PID controller keeps $A'$ stable as well, the system will converge perfectly to $y=0$ as with $A$. This is especially interesting for practical control design, since state-dependent terms are harder to estimate with an observer than just constants, especially if one allows arbitrary forms of $F(y)$ (as we show still works). Bias introduced by systematic model errors (e.g.~due to linearization around non-exact solution) is also countered by such integral action. This section reports on the analog property for integral control on Lie groups.


\subsection{Convergence on compact Lie group}
 
For definiteness consider a simple second-order system:
\begin{eqnarray}
L_{g^{-1}}\tfrac{d}{dt}g & = & \xi^l \;\;, \quad \label{eq:2nd}
 \tfrac{d}{dt}\xi^l = b\, \xi^l + u^l + u^l_B
\end{eqnarray}
where $g$ is the position on the Lie group $G$ (think of a rotation matrix); $L_g$ is the left-invariant translation from the tangent space at group identity, that is the Lie algebra, to the tangent space at $g$; $\xi^l$ is the left-invariant velocity and belongs to the Lie algebra (think of a vector of rotation velocities, expressed in body frame); $b \, \xi^l$ is velocity damping for $b<0$, or a term making the open-loop system unstable if $b>0$; $u^l$ is a control input and $u^l_B$ a bias on this input (think of torques in the case of rotations, again expressed in body frame). The PID controller for this system, as introduced in e.g.~\cite{MaithripalaJournal,MaithripalaCDC} and \cite{OurSCL}, writes
\begin{eqnarray}
u^l & = & -k_P \text{grad}_g^l \phi - k_D\xi^l + k_I\, u^l_I \label{eq:2ndU}\\
\nonumber \tfrac{d}{dt}u^l_I & = & -k_P \text{grad}_g^l \phi - k_D\xi^l 
\end{eqnarray}
where $k_P, k_D, k_I$ are constant gains to be tuned. The function $\phi: G \mapsto \mathbb{R}$ is some potential with a minimum at the target $g_{\text{target}}$, and whose gradient is used as the equivalent of proportional action to push the system towards $g_{\text{target}}$. Finally $u^l_I$ is the integral \emph{of the other control inputs}, which belong to the Lie algebra i.e.~a vector space, unlike the error $g^{-1} g_{\text{target}}$ which belongs to the Lie group and hence cannot be integrated in a standard way.

The novelty in this section is simply that we allow $u^l_B$ to depend on $g$. We then get the following extension of our result in \cite{OurSCL}.
 
\noindent \textbf{Proposition 1:} \emph{Consider the system \eqref{eq:2nd} with PID controller \eqref{eq:2ndU}. If the bias $u^l_B(g)$ has a bounded gradient w.r.t.~$g$, then at least for $k_I>0$ and $k_D-k_I>0$ both large enough the configuration $g$ converges to the set of critical points of $\phi$, with $\xi^l=0$ and $k_I u^l_I = -u^l_B(g)$. Only the minima of $\phi$ are stable equilibria.}

\noindent \emph{Proof:} Consider the Lyapunov function candidate
$$V = \phi + \tfrac{\alpha}{2}\, \Vert \xi^l \Vert^2 + \tfrac{\beta}{2} \Vert k_I (u^l_I - \xi^l) + u^l_B \Vert^2 \, .$$
Its time derivative along trajectories writes
\begin{eqnarray*}
\tfrac{d}{dt} V & = & \text{grad}_g^l \phi \cdot \xi^l + \alpha \xi^l \cdot \dot{\xi}^l + \\
& & \beta (k_I (u^l_I - \xi^l) + u^l_B )\cdot(k_I (\dot{u}^l_I - \dot{\xi}^l) + \text{grad}_g^l u^l_B \cdot \xi^l)
\end{eqnarray*}
where we have used the notation $\dot{x} = \tfrac{d}{dt}x$. Taking $\alpha = 1/k_P$ and grouping terms, we get
\begin{eqnarray*}
\tfrac{d}{dt} V & = &  - \beta k_I\, \Vert k_I (u^l_I - \xi^l) + u^l_B \Vert^2
- \frac{(k_D-k_I-b)}{k_P}\, \Vert \xi^l \Vert^2 \\
& & + \left( k_I (u^l_I - \xi^l) + u^l_B \right) \cdot \\ & & \phantom{KKK} \left( (\tfrac{1}{k_P}-(b+k_I)\beta k_I) I + \beta \text{grad}_g^l u^l_B \right) \cdot \xi^l \, ,
\end{eqnarray*}
with $I$ the identity matrix. This is a form of the type 
\begin{equation}\label{eq:matrix}
\left[\begin{array}{l} x\\ y
\end{array}\right]^T \; \left[\begin{array}{cc}-a_1 I & A \\ A^T & -a_2 I \end{array}\right] \; \left[\begin{array}{l} x\\ y
\end{array}\right] \quad \text{ with }
\end{equation}
$a_1 = \beta k_I$, $a_2 = \frac{k_D-(b+k_I)}{k_P}$, $A = (\tfrac{1}{2k_P}-\tfrac{(b+k_I)\beta k_I}{2})\, I + \tfrac{\beta}{2} \, \text{grad}_g^l u^l_B \,$ and
$x=k_I (u^l_I - \xi^l) + u^l_B \; , \;\; y = \dot{\xi}^l \, .$
Now any condition ensuring the matrix in \eqref{eq:matrix} to be negative definite would be sufficient to conclude our proof. To show that the controller can always be tuned in this way, we here derive (loose) bounds using the Gershgorin disk theorem. Let $d_r$ (resp.~$d_c$) be a bound on the sum of absolute values on any row (resp.~column) of the matrix $\text{grad}_g^l u^l_B$.
Then for Gershgorin we need 
\begin{eqnarray}
\label{eq:cond2} k_I & > & f + d_r/2 \;\; \text{ and}\\
\nonumber k_D & > & b+k_I \; + \; k_P \, f + \beta k_P d_c/2 \;\; \\
\nonumber \text{ with} & & f = \vert \tfrac{1}{2k_P\beta}-\tfrac{(b+k_I)k_I}{2} \vert \, .
\end{eqnarray}
Take any $k_I$ satisfying $k_I>d_r$ and $k_I+b > 0$. By adjusting $\beta$, we can make $f=0$. Then the first condition of \eqref{eq:cond2} is satisfied, and the second one is easy to satisfy with $k_D$ large enough. Note that $\beta$ must be known exactly only to define the corresponding Lyapunov function, not to tune the controller gains.

With these conditions, the matrix form \eqref{eq:matrix} is negative definite. Hence $V$ is a true Lyapunov function, that stops decreasing only when $\xi^l=0$ and $k_I (u^l_I - \xi^l) + u^l_B = k_I u^l_I + u^l_B = 0$. To keep these conditions invariant, as requested by the LaSalle invariance principle, we see from \eqref{eq:2ndU} that we need $\text{grad}_g^l \phi=0$. Thus the system converges towards the set of critical points of $\phi$, with zero velocity and $u^l_B$ compensated by the integral controller. If there is a point $(\tilde{g},\tilde{u}^l_I,0)$ in the neighborhood of $(g,u^l_I,0)$ where $k_I \tilde{u}^l_I + u^l_B(\tilde{g}) = 0$ and $\phi(\tilde{g}) < \phi(g)$, then we have $V(\tilde{g},\tilde{u}^l_I,0) < \phi(g)$ and the system starting at $(\tilde{g},\tilde{u}^l_I,0)$, with decreasing $V$, can never converge back to $g$; i.e.~the equilibrium is unstable as soon as $g$ is not a minimum of $\phi$.
\hfill $\square$\\

\emph{Nota Bene:} The convergence towards the set of critical points of $\phi$ is the same as the basic result obtained without bias nor integral control. In that case already, it is unavoidable that some initial conditions, ``launched with the right speed and in the right direction'', converge towards critical points that are not minima of $\phi$; but only the minima of $\phi$ are stable. The sets of initial conditions converging to the unstable equilibria however might differ in the ``bias and integral control'' case, with respect to the nominal case.
\newline The conditions on $k_I$ and $k_D$ are only sufficient. \\

For completeness, note that for a first-order system
\begin{eqnarray}
\label{eq:1storder} L_{g^{-1}}\frac{d}{dt}g & = & -k_P \text{grad}_g^l \phi + u^l_B(g) + k_I u_l^I \\
\nonumber \tfrac{d}{dt} u^l_I & = & -k_P \text{grad}_g^l \phi
\end{eqnarray}
we have the following sufficient result.

\noindent \textbf{Proposition 2:} \emph{Denote $R$ a compact region of the state space containing a minimum $g_{\text{target}}$ of $\phi$, for which $\phi(g) = \phi_{\partial R}$ is constant for all $g \in \partial R$  the boundary of $R$, and where inside $R$ we have $\phi(g) < \phi_{\partial R}$ and $\text{Hess}_g^l\phi$ positive definite. Take $k_P$ sufficiently large such that $k_P\text{Hess}_g^l\phi - \text{grad}_g^l u^l_B$ remains positive definite inside $R$, for all expected $u^l_B$. For a set $S \subset R$, select $k_I$ large enough such that 
$$\max_{g \in S} k_I k_P \phi(g) + \tfrac{1}{2} \Vert -k_P \text{grad}_g^l\phi(g) + u^l_B(g) \Vert^2 < k_I k_P \phi_{\partial R} \, ,$$
for all expected $u^l_B$. Then the system \eqref{eq:1storder} starting at any $g(0) \in S$ will never leave $R$, and converge towards $g_{\text{target}}$.}

The proof is based on the Lyapunov function $V = k_I k_P \phi(g) + \tfrac{1}{2} \Vert -k_P \text{grad}_g^l\phi + u^l_B(g) + u^l_I \Vert^2$. A different Lyapunov function might allow to weaken the conditions of the Proposition. Note that for stability of the corresponding linear system, a similar condition on the Hessian would be \emph{necessary}, whereas the sets $R,S$ can be all $\mathbb{R}^N$.
 
 
\subsection{Example: nonlinear pendulum}\label{ssec:pendulum}

The following problem was mentioned to us in \cite{RodTutorial} as a stumbling point for standard passivity-based stabilization. Consider the pendulum in a vertical plane whose motion is governed by
\begin{equation}\label{eq:penddyn}
\tfrac{d^2}{dt^2} \theta = b \tfrac{d}{dt}\theta - w\, \sin(\theta) + u
\end{equation}
where $\theta$ is the angular position of the pendulum, with $\theta=0$ pointing downwards; $w$ is the ratio between gravitational acceleration and moment of inertia of the pendulum; $b<0$ represents some friction; and $u$ is the control torque, in appropriately normalized units. The task is to \emph{robustly} stabilize the pendulum at some angle $\theta_{\text{target}} \neq 0$, despite the gravitational potential attracting it towards $\theta=0$. Due to the latter, proportional-derivative action is not enough. Of course with perfect knowledge of the parameters one can easily reshape the gravitational potential to have a minimum at $\theta_{\text{target}}$, but any error in parameters will then imply a mismatch in $\theta_{\text{target}}$. We next show that with PID control, stabilization perfectly at $\theta_{\text{target}}$ can be achieved even with only approximate knowledge of the system parameters --- in fact we only use bounds on $w$ and $b$.

Due to the fact that after moving by $2\pi$ one is back at the original location, the circle is a nonlinear space, and in fact a Lie group. This can be highlighted by defining $g=e^{i\theta}$, such that indeed $g_1 g_2 = e^{i(\theta_1+\theta_2)}$ is another point on the circle, $\theta=0$ corresponds to the identity for multiplication, and $g^{-1} = e^{-i\theta}$ exists for all $g$ and satisfies all smoothness requirements. We can then further define $\omega = \tfrac{d}{dt}\theta$ and rewrite \eqref{eq:penddyn} as: 
\begin{eqnarray*}
e^{-i \theta} \tfrac{d}{dt} e^{i \theta} = \omega \;\;, \quad
\tfrac{d}{dt} \omega = b \omega + u + u_B
\end{eqnarray*}
where we view $u_B = - w\, \sin(\theta)$ as an undesired state-dependent bias. Here we have just dropped the superscript $l$ and replaced $\xi$ by $\omega$ compared to the general notation. To stabilize $\theta_{\text{target}}$ we can select the potential $\phi = \sin(\frac{\theta-\theta_{\text{target}}}{2})^2$ whose only extrema are a minimum at $\theta_{\text{target}}$ and a maximum at $\theta_{\text{target}}+\pi$. Note that $\phi$ is equivalent to the gravitational potential, but turned as if the center of the Earth was in the direction $\theta_{\text{target}}$. Proposition 1 then readily ensures that PID control, with appropriate tuning, will make the system converge to either $\theta=\theta_{\text{target}}$ or $\theta=\theta_{\text{target}}+\pi$, of which only $\theta_{\text{target}}$ is stable. The tuning requirements \eqref{eq:cond2} only require the bound $w$ on $\tfrac{d}{d\theta} u_B$, and possibly a bound on $b$.

\begin{figure*}
\includegraphics[width=120mm]{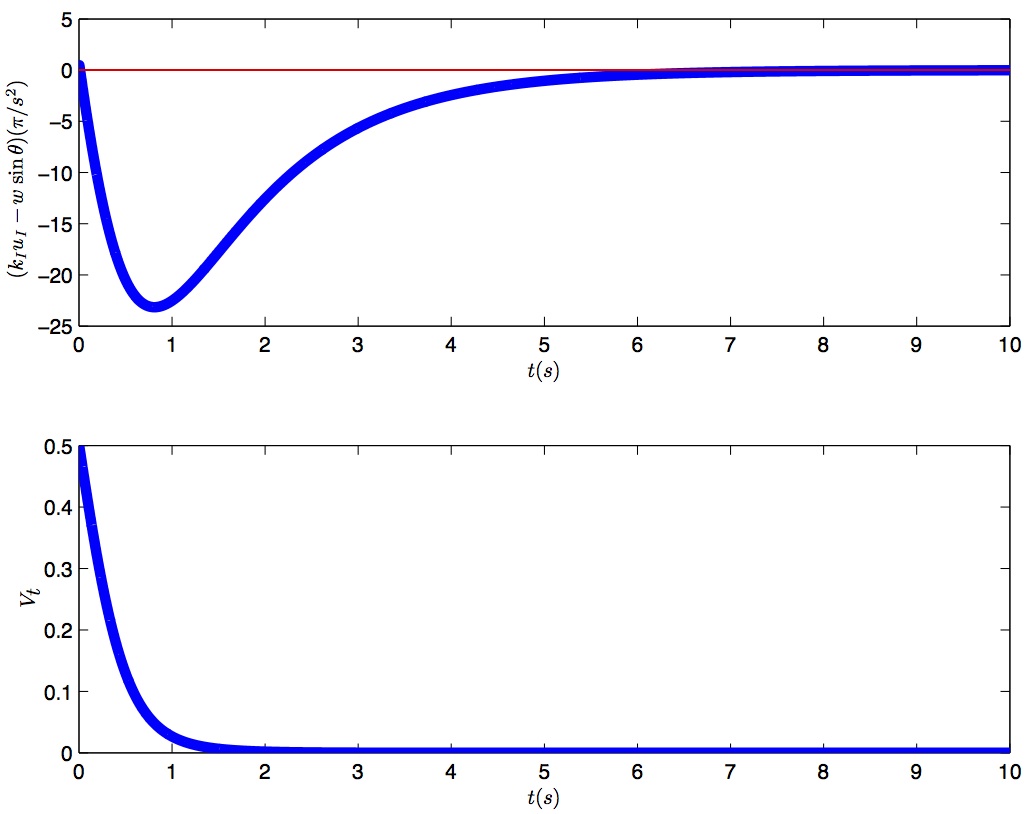} \hfill \includegraphics[width=120mm]{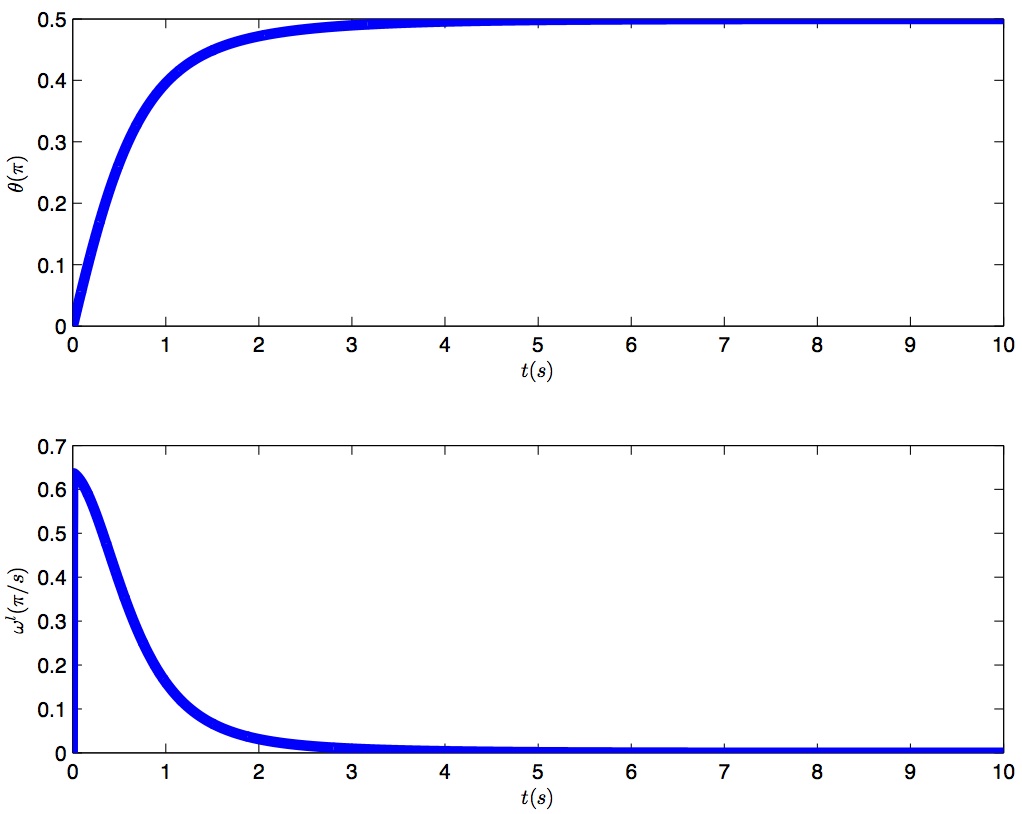}
\caption{Perfect stabilization by PID control of $\theta_{ref}=\pi/2$, despite the action of the gravitational potential considered as a priori unknown. The parameters are $w=1$, and we take $b=100$, $k_P=1000$, $k_I=1$, $k_D=600$ to satisfy \eqref{eq:cond2}. The Lyapunov function is shown with $\alpha=1/k_P=0.001$ and $\beta=1/(k_P*(b+k_I)*k_I) = 1/101000$.\vspace{4mm}}
\label{fig1}
\end{figure*}

Figure \ref{fig1} shows a simulation of this example with initial conditions $\theta(0)=\omega(0)=u_I=0$, perfectly stabilizing the state $\theta_{ref}=\pi/2$. Admittedly, the circle is not the most challenging Lie group. However, we believe that this example illustrates how  also on nonlinear spaces, integral control can be a standard method to solve stabilization problems.

\section{The steering-controlled vehicle, 1:\newline without integral action}\label{sec:car}

Consider a planar rigid body, with position $p \in \mathbb{R}^2$ and orientation given by a $2\times 2$ rotation matrix $Q$. The dynamics for its \emph{steering control} writes \cite{Krishnaprasad,RodolphePaley}:
\begin{equation}\label{Model}
\tfrac{d}{dt}Q\;=\; Q\,\cdot Q_{\frac{\pi}{2}}\,\cdot \omega \;\; , \quad \tfrac{d}{dt}p\;=\; Q\, v
\end{equation}
where $Q_{\frac{\pi}{2}}$ is a $\pi/2$ rotation in the plane; $\omega \in \mathbb{R}$ is the angular velocity, \emph{to be controlled}; and $v \in \mathbb{R}^2$ is the translation velocity expressed in body frame, \emph{which is fixed} both in direction and in magnitude. The fixed direction of $v$ can be interpreted as a vehicle which cannot translate sidewards. The fixed magnitude can be for technological simplification (on-off motor, propulsion system) or follow from floating constraints (airplane, buoyancy-driven underwater vehicle).

In motion tracking of the non-holonomic vehicle, references $p_r(t)$ and $Q_r(t)$ are given which satisfy \eqref{Model} for some $\omega_r(t)$ and the nominal $v$, which we denote $v_r=\mathbf{e}_1=[1,0]$ without loss of generality. For simplicity we consider $\omega_r(t)=\omega_0 >0$ constant, in which case the reference is moving at rate $\omega_0$ on a circle of radius $v/\omega_0$.

\noindent \textbf{Definition 1:} \emph{We say that the vehicle follows the reference \emph{motion} if and only if, in a frame that moves and turns with the reference, the position and orientation of the controlled vehicle are constant. In other words, the reference and the vehicle together form a virtual rigid body, of unimposed shape. For a reference with constant $\omega_r=\omega_0$ and a vehicle model perfectly matching the reference one, this means that the vehicle follows the same circular path as the reference, with a constant phase delay.}

We call this ``motion tracking'' to distinguish it from the traditional tracking problem, where the vehicle would have to reach the exact same position at time $t$ as the reference at time $t$. The advantage of just motion tracking is that in fact the unimposed shape can be controlled independently by another controller \cite{Krishnaprasad,RodolphePaley,AlainLieGroup}. Practical motivation for this is the possibility to independently optimize the shape of the virtual rigid body according to objectives like formation flying/navigation/platooning, or the carrying of a heavy rigid load by a team of vehicles.

To steer the vehicle towards motion tracking and stabilize it there, we can follow the geometric control design of 
\cite{RodolphePaley,AlainLieGroup}. Indeed, applying this motion coordination approach to the particular case of a network of two agents we get $\omega = \omega_0 + \omega_P$ with
\begin{equation}\label{PositionCorrection}
\omega_P = k_P \mathbf{e}_1^T \left(\omega_0 Q^T (p-p_r)-Q^T Q_r Q_{\frac{\pi}{2}}\mathbf{e}_1 \right)
\end{equation}
where $k_P>0$ is a proportional gain. Note that this controller follows as the vehicle measures $Q^T (p-p_r)$ its relative position to the reference, expressed in body frame, and $Q^T Q_r$ its relative orientation with respect to the reference. If the vehicle perfectly fits the nominal model, with $v =\mathbf{e}_1$, then the results of \cite{RodolphePaley,AlainLieGroup} imply that it will almost globally converge towards the reference motion. Our goal here is to ensure similar tracking performance when $v\neq \mathbf{e}_1$ is not perfectly known.

A central object in our analysis, as in \cite{RodolphePaley,AlainLieGroup}, is the point  
$$c(t)=p(t)+\frac{Q_{\frac{\pi}{2}}}{\omega_0}Q(t)v \; .$$ 
This is the center of the circle around which the vehicle would rotate if it just applied $\omega=\omega_0$ from its current position and orientation. Note that the exact $c$ is in fact unknown to the vehicle if it does not know its $v$ exactly. The convergence proof in \cite{RodolphePaley,AlainLieGroup} is based on the distance between $c(t)$ and the center $c_r$ of the reference circle. The combination of $c(t)$ converging to $c_r$ and $\omega$ converging to $\omega_0$, as the need for position correction $\omega_P$ asymptotically vanishes, eventually implies that the final behavior satisfies motion tracking.


\subsection{Magnitude bias}

When $v= (1+\epsilon) \mathbf{e}_1$ differs from $v_r$ but is still parallel to it, the system with just the nominal controller \eqref{PositionCorrection} in fact behaves very well. The proof is based on the same Lyapunov function as for the nominal case.\vspace{2mm}

\noindent \textbf{Proposition 3:} \emph{Consider the system \eqref{Model} with $v = (1+\epsilon) \mathbf{e}_1$ and steering control $\omega = \omega_0 + \omega_P$, where $\omega_P$ is defined by \eqref{PositionCorrection}. Then the vehicle will converge towards motion tracking of the reference. I.e.~it will rotate at rate $\omega_0$ on a circle of center $c_r$, but of radius $v/\omega_0$ possibly different from the reference radius.}\vspace{2mm}

\noindent \emph{Proof:} Consider the candidate Lyapunov function 
\begin{eqnarray}\label{eq:Vit}
		V(Q,p) &=&\tfrac{1}{2}\parallel c(t)-c_r\parallel^2\\
\nonumber	&=&\tfrac{1}{2}\parallel p(t)-p_r(t)+\tfrac{Q_{\frac{\pi}{2}}}{\omega_0}(Q(t)v-Q_r(t)\mathbf{e}_1)\parallel^2 \, .
\end{eqnarray}
Since $c_r$ is constant, when computing $\frac{d}{dt}V$ one only needs
\begin{equation}\label{eq:dV1}
	\begin{split}
		\tfrac{d}{dt}c\;=\;& \tfrac{dp}{dt}+\tfrac{Q_{\frac{\pi}{2}}}{\omega_0}\tfrac{dQ}{dt} v
		\;=\; Qv + \tfrac{Q_{\frac{\pi}{2}}}{\omega_0}Q\,Q_{\frac{\pi}{2}}(\omega_0+\omega_P)v\\
		\;=\;&\frac{-\omega_P}{\omega_0}\, Q v \, ,
	\end{split}
\end{equation}
where we have used the fact that all planar rotations commute and $Q_{\frac{\pi}{2}} Q_{\frac{\pi}{2}} v = Q_{\pi} v = -v$. Rewriting $\omega_P$ in terms of $c-c_r$ instead of $p-p_r$, see \eqref{eq:Vit}, we get
\begin{equation}\label{eq:omegapinc}
		\omega_P = k_P \omega_0\,  \left(\mathbf{e}_1^T Q^T(c-c_r) \right)\, - k_P \mathbf{e}_1^T Q_{\frac{\pi}{2}} v.
\end{equation}
The last term drops when $v$ is parallel to $\mathbf{e}_1$. Replacing this in \eqref{eq:dV1} yields
\begin{equation}
	\begin{split}
		\tfrac{d}{dt}V\;=\;&(c-c_r)^T\tfrac{d}{dt}(c-c_r) \;=\;(c-c_r)^T\tfrac{d}{dt}c\\
		\;=\;& -k_P \left(\mathbf{e}_1^T Q^T(c-c_r) \right) \; (c-c_r)^T \, Qv \\
		\;=\;& -k_P |v| \left(\mathbf{e}_1^T Q^T(c-c_r) \right)^2
		\;\leq\; 0 \, .
	\end{split}
\end{equation}
Now let us apply the LaSalle invariance principle. The set characterized by $\frac{d}{dt}V=0$ is $S=\{(p,Q)|\omega_P=0\}=\{(p,Q)|c=c_r \,\text{or}\, \cos\theta=0\}$, where $\theta$ is the angle between the vectors $Q^T(c-c_r)$ and $\mathbf{e}_1$. The subset $S_1=\{(p,Q)|c=c_r\}$ is invariant since $\omega_P=0$ implies $dc/dt=0$ from \eqref{eq:dV1}. On the other hand, starting from any point in the subset $S_2=\{(p,Q)|c\neq c_r\,\text{and}\, \cos\theta=0\}$, we will have 
$$\tfrac{d}{dt}\cos\theta = - \omega_0 \mathbf{e}_1^T Q_{\frac{\pi}{2}} Q^T(c-c_r) \; .$$
 But $Q_{\frac{\pi}{2}} Q^T(c-c_r)$ belongs to the same plane as $Q^T(c-c_r)$ and $\mathbf{e}_1$, and it is orthogonal to $Q^T(c-c_r)$. Thus when $\mathbf{e}_1^T Q^T(c-c_r) = 0$ and $(c-c_r)\neq0$ we have $\mathbf{e}_1^T Q_{\frac{\pi}{2}} Q^T(c-c_r) \neq 0$, such that a situation where $\cos\theta=0$ cannot be invariant. This concludes the proof. \hfill $\square$\\

Proposition 3 thus implies that, even if the actual translation velocities of different vehicles are unknown, then just applying the same proportional feedback \eqref{PositionCorrection} ensures that they will naturally converge towards circles of the same center, and of radii exactly matched such that they move as a rigid formation.


\subsection{Velocity misalignment: effect with nominal controller}

We now investigate what happens under nominal control when $v$ is not parallel to $\mathbf{e}_1$, i.e.~there is some misalignment by a rotation $Q_\phi$ in the constant propulsion direction of the vehicle. The ideal control would then simply apply \eqref{PositionCorrection} after correcting the misalignment in the definition of the body frame of the vehicle, i.e.~
\begin{eqnarray*}
\omega^\text{ideal}_P & = & k_P \mathbf{e}_1^T Q^T Q_\phi^T \left(\omega_0 (p-p_r)-Q_r Q_{\frac{\pi}{2}}\mathbf{e}_1 \right) \\
& = & k_P \omega_0 v^T Q^T (c-c_r) \, .
\end{eqnarray*}
From the previous section, a possible remaining difference in velocity magnitude has no detrimental effect. The difference
\begin{eqnarray*}
\omega_B & = & \omega_P-\omega^\text{ideal}_P \\
& = & k_P \mathbf{e}_1^T Q^T (I-Q_\phi^T) \left(\omega_0 (p-p_r)-Q_r Q_{\frac{\pi}{2}}\mathbf{e}_1 \right)
\end{eqnarray*}
can be viewed as a state-dependent actuation bias.
%

\begin{figure*}
\includegraphics[width=120mm]{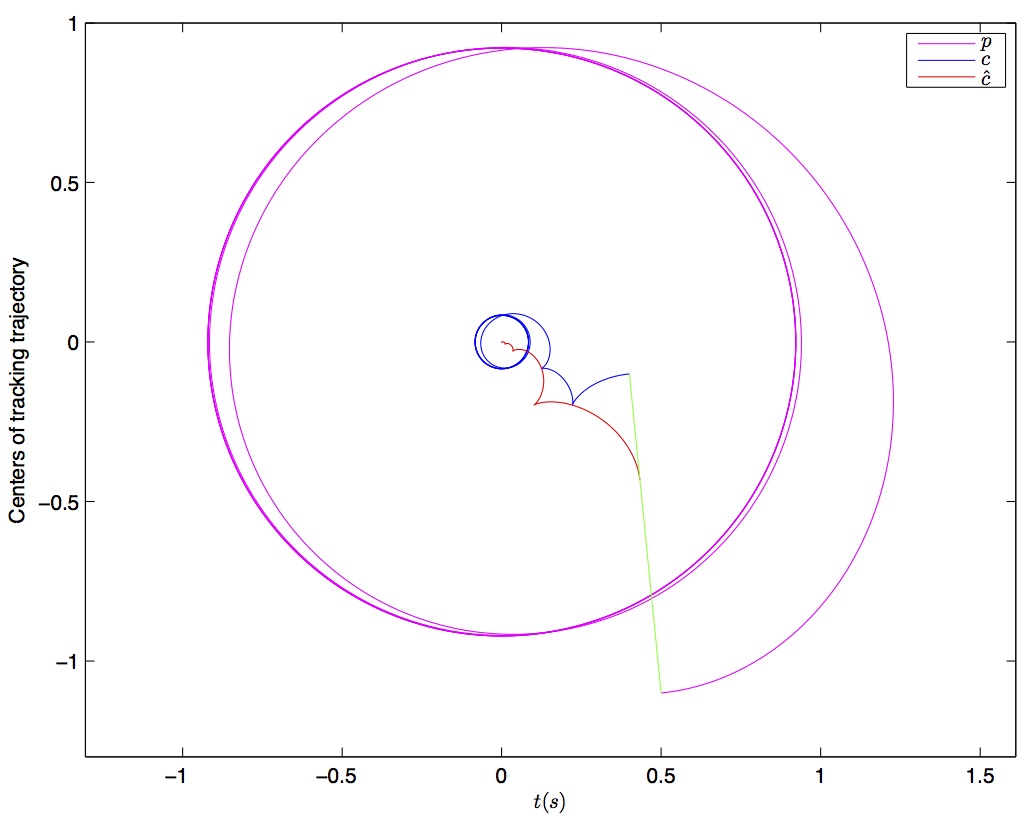} \hfill \includegraphics[width=120mm]{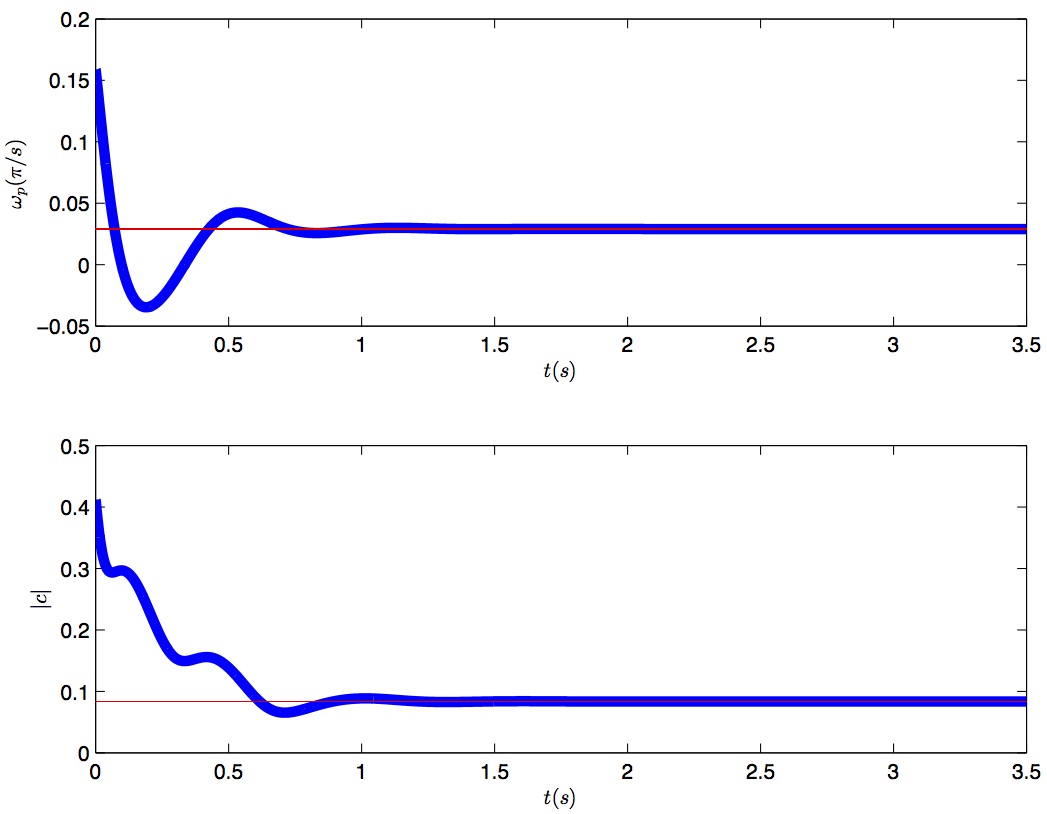}
\caption{Vehicle motion under steering control, with nominal proportional-like controller $\omega=\omega_0+\omega_P$ (taking $\omega_0=k_P=1$) and with translation velocity misalignment $v=[1,0.1]$. The left plot shows the trajectories of position $p(t)$, circle center $c(t)$ computed as if rotating with $\omega=\omega_0$ (ideal final value), and circle center $\hat{c}(t)$ computed as if rotating with $\omega=\omega_0+\hat\omega_P$ (actual final value); the initial values are marked by the green line. The right plot shows $\omega_p$ converging towards $\hat\omega_P = -\omega_0/2 + \sqrt{\omega_0^2/4 + \omega_0 k_P \cdot 0.1} =0.0916$ instead of towards $0$, and the corresponding distance between $c(t)$ and the target center $c_r=0$.}\label{fig:2}
\end{figure*}

Simulations, see Fig.\ref{fig:2}, show that when $\omega_P$ is applied while $\phi \neq 0$, the center $c(t)$ converges to a circular limit cycle around $c_r$. Moreover, $\omega$ converges towards a constant value $\hat\omega = \omega_0 + \hat\omega_P$ with $\hat\omega_P \neq 0$. 
This is another viewpoint about $\hat\omega_P$ as an input bias on $\omega_0$.

Although a proof of convergence towards this situation is beyond the scope of the present paper, we can easily compute the value of $\hat\omega_P$ by examining the conditions for the equilibrium observed in simulations. 
Note that the center around which the vehicle is turning is given by $\hat{c} = p(t)+\frac{Q_{\frac{\pi}{2}}}{\omega_0 + \hat\omega_P}Q(t)v$. Constant rotation on a circle of fixed center $\hat{c}$ requires $\tfrac{d}{dt}\omega_P = 0$ and $d\hat{c}/dt=0$. Expressing $\omega_P$ in terms of $\hat{c}$ and $c_r$ instead of $p,p_r$ and imposing $\hat{c}=c_r$, we get the condition
$$\omega_P = \tfrac{\omega_0}{\omega_0+\omega_P} k_P  \vert v \vert \sin\phi \, .$$
Solving this, we get
\begin{equation}
\hat\omega_P = \tfrac{-\omega_0}{2} \pm \sqrt{\tfrac{\omega_0^2}{4} + \omega_0 k_P \vert v \vert \sin\phi} \; .
\end{equation}
One checks (see captions) that the solution with the plus sign fits the value observed in simulations.

This is a somewhat strange and undesirable ``steady state'', as the rotation speed is not the reference one and thus $Q^T Q_r$ is time varying, as well as the position of the vehicle with respect to the reference! Then the reference and following vehicle will not move at all like a rigid body. Instead, in a frame attached to the reference, the follower will move on a circle of radius $\bar{v}/{\bar{\omega}}$, thus periodically coming possibly very close and getting very far from the reference. In other words, the result of the bias is a small, constant but thus accumulating drift. This might be easily counteracted by the formation shape controller. Our goal now however is to counteract this effect directly, without imposing a shape to the formation, by introducing appropriate integral action.

 
\section{The steering-controlled vehicle, 2: adapted integral action}\label{sec:car2}

We now consider the integral action to counteract the velocity misalignment. The difference $\omega_P - \omega^\text{ideal}_P$ induced by $\phi\neq 0$ is equivalent to a state-dependent actuation bias, i.e.~we actually \underline{apply} some $u(x)+u_B(x)$ (here $\omega_P$) while we \underline{want} to apply some $u(x)$ (here $\omega^\text{ideal}_P$). However, towards applying additional corrective actions, there is an important difference: in an input bias situation we \underline{know} $u(x)$, while here what we know is $u(x)+u_B(x)$. This corresponds more to the effect of an \emph{output bias}, e.g.~in a linear system:
\begin{eqnarray*}
\tfrac{d}{dt}x & = & Ax + B u\;\; , \;\; y = (C+C_B) x \;\; , \;\; u=-K\,y\\
\Rightarrow & &  \tfrac{d}{dt}x = (A-BKC) x - BKC_B x \, .
\end{eqnarray*}
The presence of $C_B$ has an effect equivalent to a state-dependent input bias $KC_B x$; but what we can use for further feedback control, thus $y$, is related to $(C+C_B)\, x$, thus the command \emph{including} the bias.

It turns out that an integral-like action \emph{does} help reject the effect of this output bias, at least in the present application (where it is equivalent to an actuation misalignment). 


\subsection{Integral action design and convergence analysis}\label{sec:}

We make two changes in the integral control.
\begin{itemize}
\item[1.] We reverse the sign of integral control, such that the integral action does not reinforce but rather damp the effect of the actuation. This can be understood by remembering indeed that this actuation \emph{includes the bias}.
\item[2.] We integrate all input commands, including the integral correction itself:
\begin{equation}\label{eq:intoo}
\begin{split}
\omega \; = \;& \omega_0 + \omega_P - k_I \omega_I \;\; \text{with}\\
\tfrac{d}{dt} \omega_I \; =\;& \omega_P - k_I\omega_I\\
\; =\;& \omega^\text{ideal}_P + \omega_B - k_I\omega_I\;.
\end{split}
\end{equation}
\end{itemize}

Note that with this notation, as soon as $\tfrac{d}{dt}\omega_I=0$, we have the desired situation $\omega = \omega_0$; all such situations are not necessarily invariant, and this is fortunate. Indeed, to completely reach our target, we must also get to the correct circle center. The situation $c=c_r$ (thus $\omega^\text{ideal}_P=0$) and $k_I\omega_I = \omega_B$ is invariant. We can give the following convergence property towards it.

\noindent \textbf{Proposition 4:} \emph{For any $\vert v \vert$ and any misalignment $\phi \in (-\tfrac{\pi}{2},\tfrac{\pi}{2})$, there exist $k_I$ and $k_P$ small enough such that the target situation, with $\omega=\omega_0$ and $c=c_r$, is locally asymptotically stable with controller \eqref{eq:intoo}.}

\noindent \emph{Proof:} The system becomes simpler in a frame rotating at $\omega$. Hence we define $x = Q^T\, c$
and we get
\begin{eqnarray*}
\tfrac{d}{dt} x & = & \tfrac{-v}{\omega_0} \left[k_P\omega_0 \mathbf{e}_1^T x - k_P \mathbf{e}_1^T Q_{\frac{\pi}{2}} v - k_I \omega_I \right] - \omega Q_{\frac{\pi}{2}} x\\
\tfrac{d}{dt} \omega_I & = & \left[k_P\omega_0 \mathbf{e}_1^T x - k_P \mathbf{e}_1^T Q_{\frac{\pi}{2}} v - k_I \omega_I \right] \, .
\end{eqnarray*}
This system features an equilibrium at $x=0$, $\omega_I = \tfrac{-k_P}{k_I} \mathbf{e}_1^T Q_{\frac{\pi}{2}} v = \tfrac{k_P}{k_I} \vert v \vert \sin\phi =: \hat\omega_I$. It is nonlinear only through the term $\omega Q_{\frac{\pi}{2}} x$, which close to the equilibrium linearizes to $\omega_0 Q_{\frac{\pi}{2}} x$.
The linearized system then has the characteristic polynomial
\begin{eqnarray*}
P(\lambda)& = & \lambda^3 + (k_I + \vert v \vert k_P \cos\phi) \lambda^2 \\
& & + (\omega_0^2 + \omega_0 \vert v \vert k_P \sin\phi) \, \lambda + k_I \omega_0^2 \, ,
\end{eqnarray*}
which according to the Routh-Hurwitz criterion is positive for $\phi \in (-\tfrac{\pi}{2},\tfrac{\pi}{2})$ whenever
\begin{eqnarray}\label{eq:bounds1}
\omega_0 & > & k_P \vert v \vert\, \vert\sin\phi \vert \quad \text{and} \\
\nonumber \omega_0\cos\phi & > & k_I \vert \sin\phi \vert + \vert v \vert k_P \vert \sin\phi\cos\phi \vert \, .
\end{eqnarray}
For any given $\phi$ strictly inside $(\tfrac{\pi}{2},\tfrac{\pi}{2})$ the left hand sides are fixed positive constants, so it suffices to take $k_I$ and $k_P$ small enough on the right hand sides, to ensure that the linearized system is stable.\hfill $\square$\\

\begin{figure*}
\includegraphics[width=120mm]{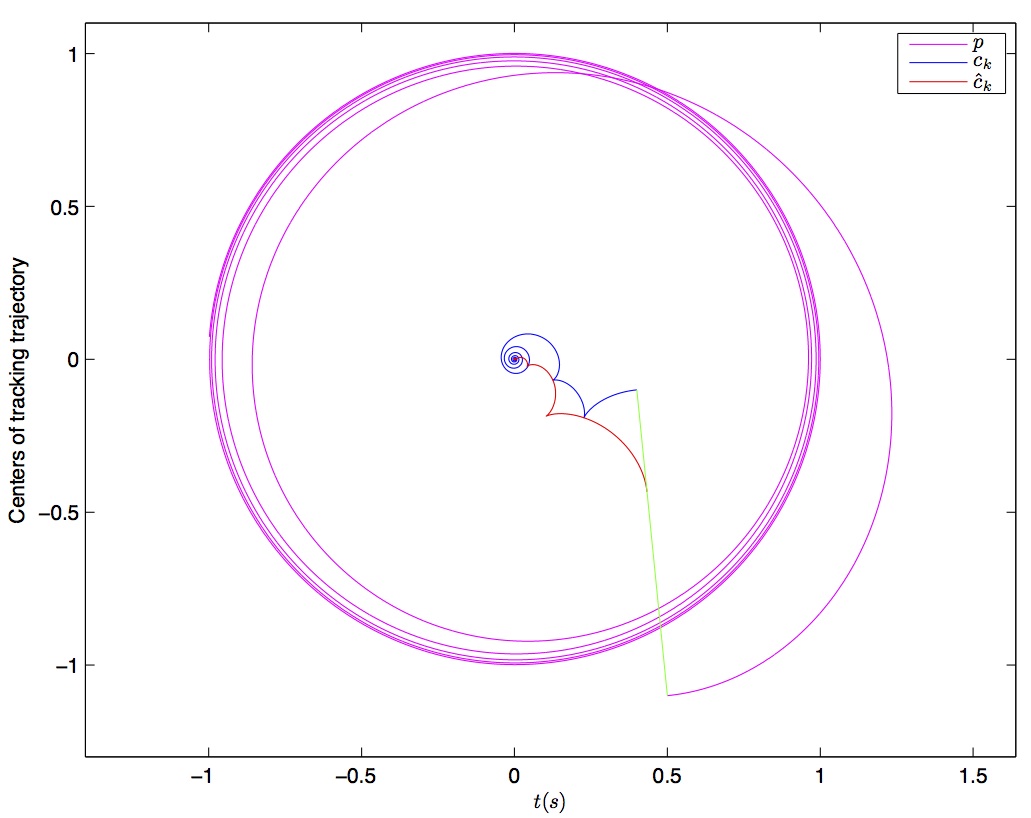} \hfill \includegraphics[width=120mm]{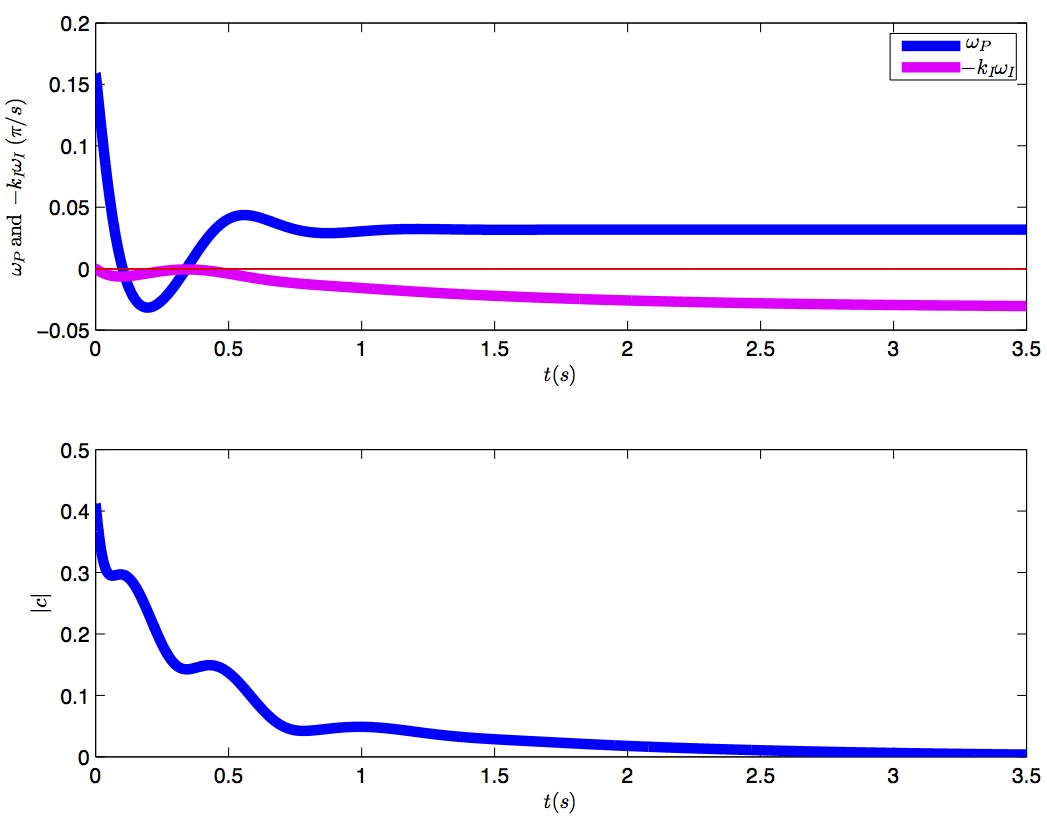}
\caption{Same plots as on Fig.\ref{fig:2}, now including the adapted integral control \eqref{eq:intoo} with $k_I=0.1$. Since $\vert v \sin\phi \vert = 0.1 = \tan\phi$, one easily checks that the conditions \eqref{eq:bounds1} are satisfied. The vehicle now converges perfectly to the target trajectory, while the bias effect is countered by the integral control input.}\label{fig3}
\end{figure*}

Figure \ref{fig3} shows a simulation with the same parameters as on Fig.~\ref{fig:2}, but adding our integral control with $k_I=0.1$. It confirms the improved behavior of the system. 

\noindent \emph{Remark:} From the presence of $\omega_0$ in \eqref{eq:bounds1}, one clearly sees the key role played by the fact that in steady state the vehicle is rotating with a constant drift input, sufficiently fast with respect to the feedback. This makes the strategy look much like an averaging scheme.




\section{Concluding remarks}

The aim of this paper is to revive the use of integral control to robustly reject biases on system dynamics. This idea has recently been formalized for constant biases on Lie groups, and we provide two extensions. First, we show that a state-dependent bias can also be perfectly rejected. This makes a strong case for using integral control in systems where natural dynamics must be cancelled, complementing the often-used open-loop model-based cancellation.  Second, we treat a case of input bias on an underactuated system, namely steering control of a vehicle moving at constant speed. When the bias is on the speed magnitude, it seems to have no detrimental effect on keeping the vehicle in a formation with fellows. When the bias is on the direction however, it has an effect and must be compensated. We had to modify the integral controller to cope with this case, which in fact appears closer to a bias in the measurement output than in the control input.

This last point suggests, for future research, to explore further links between systems with output bias and systems with underactuated input bias, to possibly solve one case with tools from the other. E.g.~one might want to investigate in which cases integral control can help to efficiently and simply reject an output bias. The conditions for local asymptotic stability (Proposition 4) are analog to those of averaging schemes, and in fact there is hope to prove global convergence using a combination of the exact local stability proved here, with an approximate but global trajectory characterization using the averaging theorems of nonlinear dynamical systems. This averaging behavior also gives indications as to how the integral controller is ``averaging out'' the bias, a viewpoint which might be useful for the design of other bias-rejecting controllers on nonlinear spaces. 

To conclude, we must note that velocity bias cannot always be exactly compensated, see e.g.~our benchmark model with $\omega_0=0$ and a wrong magnitude of $v$; hence the full treatment of bias rejection for underactuated systems on Lie groups remains an open question to our knowledge.

\section*{Acknowledgments}

This paper presents research results of the Belgian Network DYSCO  
(Dynamical Systems, Control, and Optimization), funded by the Interuniversity Attraction Poles Programme, initiated by the Belgian State, Science Policy Office. The first author's visit to Ghent University has been supported by a CSC scholarship, initiated by the China Scholarship Council. The authors want to thank prof.R.Sepulchre for suggesting the pendulum example.

\end{document}